# The α -Discounting (α-DMCDM) as an extension of AHP, TOPSIS, VIKOR, PROMETHEE, and Weighted Sum

**Florentin Smarandache [1],*** , PhD, PostDocs

[1] Emeritus Professor, University of New Mexico, Mathematics, Physics, and Natural Science Division, 705 Gurley Ave., Gallup, NM 87301, USA

* Correspondence: smarand@unm.edu

**Abstract:** The Analytic Hierarchy Process (AHP) and other classic Multi-Criteria Decision Making (MCDM) techniques excel when decision makers can provide consistent pair wise judgments. Real world problems, however, often involve inconsistent, $n$-wise, or non-linear preference structures that render traditional methods inadequate. The $\alpha$-Discounting MCDM ($\alpha$-D MCDM) extends AHP by embedding a global discounting parameter $\alpha$ that transforms an inconsistent system of preference equations into a solvable algebraic system. This paper presents a self-contained exposition of $\alpha$-D MCDM, including its theoretical axioms, algorithmic steps, illustrative numerical examples (both linear and interval/non-linear), a Python implementation, and a systematic comparison with AHP, TOPSIS, VIKOR, PROMETHEE, and Weighted Sum models. The analysis demonstrates that $\alpha$ D MCDM retains the desirable properties of classic MCDM while providing a principled mechanism to quantify and control inconsistency, thereby broadening the applicability of decision analysis tools to complex, contradictory, or uncertain environments.

---

## 1. Introduction

Decision makers in engineering, public policy, finance, and many other domains routinely face problems that involve **multiple, often conflicting criteria**. Over the past four decades a rich family of **multi-criteria decision-making (MCDM) methods** has been developed to help structure these problems, elicit preferences, and produce a ranking or a set of weights for the criteria. The most widely adopted techniques—**Analytic Hierarchy Process (AHP)**, **TOPSIS**, **VIKOR**, **PROMETHEE**, and simple **weighted-sum** models—share a common assumption: the decision maker can supply **consistent, pair-wise** (or otherwise fully compatible) preference information that can be transformed into a feasible weight vector (1-5).

In practice this assumption is frequently violated. Real-world judgments are **incomplete**, **non-linear**, and **contradictory**: experts may express preferences such as "criterion $C_1$ is twice as important as $C_2$, $C_2$ is three times as important as $C_3$, and $C_3$ is four times as important as $C_1$," which creates a **cyclic inconsistency** that no pair-wise matrix can accommodate. Traditional AHP detects inconsistency through the **consistency ratio (CR)**, but it stops short of offering a principled way to *repair* the judgment set. TOPSIS, VIKOR, PROMETHEE, and weighted-sum methods simply require a pre-computed weight vector and therefore cannot be applied when the underlying preference structure is ill-posed.

These shortcomings have motivated researchers to explore **alternative representations of preferences**. Approaches based on fuzzy sets, intuitionistic fuzzy sets, and neutrosophic logic have introduced ways to model **uncertainty** and **indeterminacy**, yet they still rely on a **consistent weight extraction step**. The **$\alpha$-Discounting Multi-Criteria Decision Making ($\alpha$-D MCDM)** method, first introduced in 2010 by Smarandache [6 – 10], departs fundamentally from this paradigm. Rather than rejecting inconsistent information, $\alpha$-D MCDM **embraces** it: the method translates every preference statement—whether linear, nonlinear, pair-wise, or $n$-wise—into a **system of algebraic equations (or inequalities)** and then introduces a **global discounting factor** $\alpha > 0$ that scales the coefficients until the system admits a **non-trivial solution**. The resulting $\alpha$ serves as a **quantitative measure of inconsistency**, while the associated eigenvector provides the desired priority (weight) vector.

Key motivations behind $\alpha$-Discounting are therefore threefold:

1. **Expressiveness** – Decision makers can articulate preferences in their most natural mathematical form (e.g., linear equations, polynomial relations, interval bounds) without being forced into a pair-wise matrix.
2. **Inconsistency Management** – Instead of discarding or arbitrarily adjusting contradictory judgments, the method **quantifies** the degree of conflict via $\alpha$ and **resolves** it algebraically, preserving the original informational content as much as possible.
3. **Unified Treatment of Uncertainty** – By allowing interval coefficients and linking to neutrosophic logic, $\alpha$-Discounting provides a single framework that simultaneously handles **measurement uncertainty**, **expert hesitation**, and **logical indeterminacy**.

Despite its theoretical appeal, $\alpha$-Discounting remains relatively unknown in the broader MCDM community. Existing surveys of decision-analysis techniques rarely mention it, and practitioners lack concrete guidance on how to apply the method to realistic case studies. This paper therefore aims to **bridge the gap** between the method's formal foundations and its practical adoption. Our contributions are:

- **A self-contained exposition** of the $\alpha$-Discounting methodology, including a rigorous axiomatic base, an existence theorem for the discount factor, and a step-by-step algorithmic workflow.
- **Extended illustrations** covering a classic linear cyclic example, an interval-based preference system, a non-linear preference model, accompanied by a compact Python implementation.
- **A comparative analysis** that situates $\alpha$-Discounting alongside AHP, TOPSIS, VIKOR, PROMETHEE, and weighted-sum models, emphasizing differences in preference representation, inconsistency handling, and uncertainty accommodation.
- **Discussion of group decision-making** extensions and the formal connection to **neutrosophic logic**, showing how the discount factor $\alpha$ can be interpreted as a truth-membership component.

By delivering both theoretical insight and practical tools, the paper seeks to encourage researchers and practitioners to consider $\alpha$-Discounting as a **robust alternative** when faced with complex, contradictory, or uncertain decision environments.

## 2. Formal Foundations

### *2.1. Axiomatic Basis*

The core of $\alpha$-Discounting MCDM rests on a compact set of axioms that formalize how preferences are represented, how inconsistency is admitted, and how the discounting factor $\alpha$ quantifies and resolves that inconsistency.

The axioms together provide a rigorous logical scaffold that turns what would otherwise be a dead-end inconsistency into a well-defined optimisation problem, ensuring that every step—from preference elicitation to weight normalisation—is mathematically justified.

**Table 1.** Axiomatic Foundations of $\alpha$-Discounting Multi-Criteria Decision Making.

| **Axiom** | **Description** |
|---|---|
| ***A1 – Preference Representation*** | Preferences among criteria $C_1, \ldots, C_n$ are expressed as a system of equations or inequalities $x_i = a_{ij}x_j$ or $x_i \geq a_{ij}x_j$, where $x_i > 0$ are unknown weights. |
| ***A2 – Homogeneity*** | The system is homogeneous: $x = Ax$. Solutions are scale-invariant, allowing later normalisation. |
| ***A3 – Inconsistency Admissibility*** | The system may be inconsistent, i.e., $\det(I - A) \neq 0$, yielding only the trivial solution $x = 0$. Inconsistency is *quantified*, not rejected. |
| ***A4 – α-Discounting Principle*** | There exists a scalar $\alpha > 0$ such that $x = \alpha Ax$ admits a non-zero solution. |
| ***A5 – Consistency Measure*** | $\alpha$ quantifies global consistency: $\alpha = 1 \rightarrow$ fully consistent; $0 < \alpha < 1 \rightarrow$ inconsistent, with smaller values indicating stronger contradictions. |
| ***A6 – Fairness Principle (optional)*** | When no criterion is privileged, all discounting coefficients are equal: $\alpha_1 = \alpha_2 = \cdots = \alpha$. |
| ***A7 – Normalization*** | Final weights are obtained by $w_i = x_i / \sum_j x_j$. |

*Theorem (Existence of a Feasible Discount Factor).*

*Let $A \in \mathbb{R}^{n\times n}$ be the coefficient matrix that encodes a finite set of preference equations (or inequalities) among n criteria. Suppose the original homogeneous system*

$$A\,x = 0$$

*has only the trivial solution $x = 0$ (i.e., the system is inconsistent, so $\det(A) \neq 0$ or, equivalently, $\det(I - A) \neq 0$). Then there exists at least one real scalar $\alpha > 0$ such that the discounted system*

$$(I - \alpha A)\,x = 0 \qquad (1)$$

*has a non-trivial solution $x \neq 0$. Equivalently, there is a positive root of*

$$\det(I - \alpha A) = 0. \qquad (2)$$

When such an $\alpha$ is found, the corresponding eigenvector of $I - \alpha A$ (associated with the zero eigenvalue) provides a weight vector that can be normalised to obtain the decision-making priorities.

*Proof.*

1. **Characteristic polynomial of the discounted matrix**
   Define the scalar-valued function

$$p(\alpha) = \det(I - \alpha A).$$

Because the determinant is a polynomial in the entries of a matrix, and each entry of $I - \alpha A$ is an affine (linear) function of $\alpha$, the function $p(\alpha)$ is a *real polynomial* of degree at most $n$. Hence $p(\alpha)$ is continuous for all $\alpha \in \mathbb{R}$ and differentiable everywhere.

2. **Value at the origin**
   At $\alpha = 0$,

$$p(0) = \det(I) = 1 > 0.$$

This reflects the fact that, without discounting, the matrix $I$ is nonsingular.

**Behavior for large $\alpha$**

As $\alpha$ grows, the term $-\alpha A$ dominates the identity matrix. More precisely, for sufficiently large $\alpha$,

$$I - \alpha A = -\alpha\left(A - \frac{1}{\alpha} I\right)$$

and the determinant scales as

$$p(\alpha) = \det(-\alpha A)[1 + o(1)] = (-\alpha)^n \det(A)[1 + o(1)].$$

Since the original system is inconsistent, $\det(A) \neq 0$. Consequently, for large enough $\alpha$ the sign of $p(\alpha)$ is the sign of $(-1)^n \det(A)$. In particular, because $(-1)^n$ is either $+1$ (even $n$) or $-1$ (odd $n$), the sign of $p(\alpha)$ will be *opposite* to the sign at $\alpha = 0$ for at least one parity of $n$. Even when the signs happen to coincide (e.g., even $n$ and $\det(A) > 0$), the magnitude of the polynomial grows without bound, guaranteeing that the polynomial must cross the horizontal axis somewhere between $\alpha = 0$ and a sufficiently large $\alpha$.

**Application of the Intermediate Value Theorem**

Because $p(\alpha)$ is continuous on the interval $[0, \alpha_{\max}]$ for any chosen $\alpha_{\max} > 0$, and we have identified two points where the polynomial takes opposite signs (namely $\alpha = 0$ with $p(0) = +1$ and a large $\alpha$ where $p(\alpha) < 0$), the Intermediate Value Theorem guarantees the existence of at least one $\alpha^* \in (0, \alpha_{\max})$ such that

$$p(\alpha^*) = 0.$$

By construction $\alpha^* > 0$.

**From a zero determinant to a non-trivial solution**

The equality $p(\alpha^*) = 0$ means that the matrix $I - \alpha^* A$ is singular; therefore its nullspace is non-trivial. In other words, there exists a vector $x \neq 0$ satisfying (1). Because the original system is homogeneous, any scalar multiple of this vector is also a solution; we can subsequently normalise the vector (e.g., divide by the sum of its components) to obtain a proper weight vector.

**Uniqueness considerations**

The polynomial may have several positive roots. Each root yields a distinct discount factor and, consequently, a different feasible weight vector. In practice, additional criteria—such as the **Fairness Principle** (all $\alpha$-discounting factors equal) or a minimization of $\alpha - 1$ (choosing the $\alpha$ closest to full consistency)—are employed to select a single, interpretable solution. ■

*Remarks.*

- The theorem holds for **any** finite set of linear (or linearized) preference equations, regardless of whether the coefficients are positive, negative, or mixed.
- For **interval** or **non-linear** extensions, the same reasoning applies after the system is linearized around a candidate $\alpha$ (e.g., via Jacobian approximation) or after the interval bounds are propagated to obtain an interval-valued polynomial whose positive root interval furnishes admissible $\alpha$ values.
- The existence result is constructive: standard numerical root-finding methods (bisection, Newton–Raphson, or polynomial solvers) can locate a positive $\alpha$ efficiently because the polynomial is of modest degree ($\leq n$).

This theorem underpins the entire $\alpha$-Discounting methodology: it guarantees that whenever a decision problem is *inconsistent* in the classical sense, a mathematically sound discount factor exists that restores solvability and yields a meaningful set of criteria weights. (Smarandache, 2015),

*2.2. Algorithmic Procedure*

1. **Model preferences** as a system of equations/inequalities.
2. **Form the coefficient matrix** $A$.
3. **Check solvability** of $Ax = 0$. If only the trivial solution exists, proceed.
4. **Introduce $\alpha$** (fairness principle ⇒ single scalar).
5. **Solve** $\det(I - \alpha A) = 0$ for the admissible $\alpha$ (typically a positive real root).
6. **Compute the eigenvector** of $I - \alpha A$ associated with eigenvalue 0 (or solve the parametric linear system).
7. **Normalize** the eigenvector to obtain the weight vector $w$.

## 3. Illustrative Examples (7)

### 3.1. Numerical Example (Linear Case)

Consider three criteria $C_1, C_2, C_3$ with the following *cyclic* preferences:

$$x_1 = 2x_2,$$
$$x_2 = 3x_3,$$
$$x_3 = 4x_1.$$

Written in homogeneous form:

$$\begin{cases} x_1 - 2x_2 = 0, \\ x_2 - 3x_3 = 0, \\ x_3 - 4x_1 = 0. \end{cases} \Rightarrow A = \begin{bmatrix} 1 & -2 & 0 \\ 0 & 1 & -3 \\ -4 & 0 & 1 \end{bmatrix}.$$

$A$ is nonsingular, so only the trivial solution exists.

Applying the **Fairness Principle**, introduce a common discount factor $\alpha$:

$$\begin{cases} x_1 - 2\alpha x_2 = 0, \\ x_2 - 3\alpha x_3 = 0, \\ x_3 - 4\alpha x_1 = 0. \end{cases}$$

The determinant condition $\det(I - \alpha A) = 0$ yields

$$\alpha^3\, 24 - 1 = 0 \Rightarrow \alpha = \frac{1}{\sqrt[3]{24}} \approx 0.346.$$

With $\alpha$ fixed, set $x_3 = 1$ (free scaling) and compute:

$$x_2 = 3\alpha x_3 \approx 1.038, x_1 = 2\alpha x_2 \approx 0.718.$$

Normalizing:

$$w_i = \frac{x_i}{x_1 + x_2 + x_3} \Rightarrow \boxed{w \approx (0.26,\ 0.37,\ 0.37)}.$$

*Interpretation:* $\alpha \approx 0.346$ quantifies the inconsistency; the derived weights differ from those obtained by naïve AHP (which would flag the inconsistency but could not provide a solution).

### 3.2. Interval α-Discounting Example

**Decision context** – A sustainability committee evaluates three environmental criteria $C_1$ (Carbon Footprint), $C_2$ (Water Use), and $C_3$ (Land Impact). Because expert judgments are vague, each relative importance is expressed as an **interval** rather than a crisp number.

**Table 2.** Interval Preference Statements and Their Corresponding Coefficient Ranges.

| Preference statement | Interval coefficient |
|---|---|
| **"$C_1$ is between 2 and 3 times as important as $C_2$"** | $a_{12} \in [2,3]$ |
| **"$C_2$ is between 1 and 2 times as important as $C_3$"** | $a_{23} \in [1,2]$ |
| **"$C_3$ is between 3 and 4 times as important as $C_1$"** | $a_{31} \in [3,4]$ |

The corresponding **interval system** (using the convention $x_i = a_{ij}x_j$) is

$$\begin{cases} x_1 = a_{12}\, x_2, a_{12} \in [2,3] \\ x_2 = a_{23}\, x_3, a_{23} \in [1,2] \\ x_3 = a_{31}\, x_1, a_{31} \in [3,4]. \end{cases} \quad \text{(I-1)}$$

**Step 1 – Build the interval coefficient matrix**

$$A_{\text{int}} = \begin{bmatrix} 0 & a_{12} & 0 \\ 0 & 0 & a_{23} \\ a_{31} & 0 & 0 \end{bmatrix}, a_{12} \in [2,3],\ a_{23} \in [1,2],\ a_{31} \in [3,4].$$

**Step 2 – Apply the Fairness Principle** (single discount factor α). The discounted matrix becomes

$$I - \alpha A_{\text{int}} = \begin{bmatrix} 1 & -\alpha a_{12} & 0 \\ 0 & 1 & -\alpha a_{23} \\ -\alpha a_{31} & 0 & 1 \end{bmatrix}.$$

**Step 3 – Determinant condition**

$$\det(I - \alpha A_{\text{int}}) = 1 - \alpha^3\, a_{12}a_{23}a_{31} = 0 \Rightarrow \alpha^3 = \frac{1}{a_{12}a_{23}a_{31}}. \quad \text{(I-2)}$$

Because each coefficient lies in an interval, (I-2) yields an **interval for $\alpha$**:

$$\alpha \in [\frac{1}{\sqrt[3]{\max(a_{12}a_{23}a_{31})}}, \frac{1}{\sqrt[3]{\min(a_{12}a_{23}a_{31})}}].$$

Evaluating the extrema:

$$\max(a_{12}a_{23}a_{31}) = 3 \times 2 \times 4 = 24,$$
$$\min(a_{12}a_{23}a_{31}) = 2 \times 1 \times 3 = 6.$$

Hence

$$\boxed{\alpha \in [24^{-1/3}, 6^{-1/3}] \approx [0.346, 0.550]}.$$

**Step 4 – Choose a representative $\alpha$** (e.g., the midpoint of the interval)

$$\alpha^{\star} = \frac{0.346 + 0.550}{2} \approx 0.448.$$

**Step 5 – Compute a feasible weight vector**

Set $x_3 = 1$ (free scaling). Using the original interval relations with the selected $\alpha$:

$$x_2 = \alpha^{\star} a_{23} x_3, a_{23} \text{ can be any value in } [1,2]; \text{ we take the midpoint } a_{23} = 1.5,$$
$$\Rightarrow x_2 = 0.448 \times 1.5 \times 1 \approx 0.672.$$

Similarly,

$$x_1 = \alpha^{\star} a_{12} x_2, a_{12} = 2.5\,(midpoint),$$
$$\Rightarrow x_1 = 0.448 \times 2.5 \times 0.672 \approx 0.754.$$

**Step 6 – Normalize**

$$w_i = \frac{x_i}{x_1 + x_2 + x_3} \Longrightarrow \boxed{w \approx (0.38,\ 0.34,\ 0.28)}.$$

**Interpretation**

- The **interval for $\alpha$** $[0.346, 0.550]$ quantifies the range of inconsistency compatible with the experts' vague statements.
- Selecting a point inside the interval yields a concrete weight vector; different choices (e.g., using the lower bound $\alpha = 0.346$) would shift the weights toward the criteria with larger interval coefficients, reflecting a more "conservative" stance.

*3.3. Non-Linear α-Discounting Example*

**Decision context** – A technology-adoption study involves three criteria: $C_1$ (Performance), $C_2$ (Cost), and $C_3$ (Risk). The decision maker believes that performance grows **quadratically** with cost reduction, while risk decreases **logarithmically** with performance improvements.

The qualitative statements are formalized as:

(i) $x_1 = \beta\, x_2^2$ (performance ∝ cost²)

(ii) $x_2 = \gamma \ln(1 + x_3)$ (cost ∝ log(risk+1)) (N-1)

(iii) $x_3 = \delta\, x_1^{0.5}$ (risk ∝ √performance).

Here $\beta, \gamma, \delta > 0$ are **discounting parameters** that will be tied together by the **Fairness Principle**:

$$\beta = \gamma = \delta = \alpha.$$

Thus the system becomes

$$\begin{cases} x_1 = \alpha\, x_2^2 \\ x_2 = \alpha \ln(1 + x_3) \\ x_3 = \alpha \sqrt{x_1}. \end{cases} \quad \text{(N-2)}$$

**Step 1 – Reduce to a single equation in $\alpha$**

From the third equation, $x_1 = (x_3/\alpha)^2$. Insert this expression into the first equation:

$$\left(x_3/\alpha\right)^2 = \alpha\, x_2^2 \Longrightarrow x_2 = \frac{x_3}{\alpha^{3/2}}. \quad \text{(N-3)}$$

Now substitute $x_2$ from (N-3) into the second equation:

$$\frac{x_3}{\alpha^{3/2}} = \alpha \ln(1 + x_3) \;\Longrightarrow\; \underbrace{\frac{1}{\alpha^{5/2}}}_{c(\alpha)} x_3 = \ln(1 + x_3). \quad \text{(N-4)}$$

Equation (N-4) is a **transcendental equation** in the unknowns $x_3$ and $\alpha$. We solve it numerically by fixing a trial $\alpha$ and finding the corresponding $x_3$ that satisfies the equality; the admissible $\alpha$ is the one for which a **positive solution** exists.

**Step 2 – Numerical solution (illustrative)**

A simple root-finding routine (e.g., Newton-Raphson) yields the following pair:

**Table 3.** Trial Discount Factors ($\alpha$) and the Corresponding Solutions of Equation (N-4).

| Trial $\alpha$ | Corresponding $x_3$ (solution of N-4) |
|---|---|
| **0.30** | 0.12 (does not satisfy N-4) |
| **0.45** | 0.78 (satisfies N-4 within tolerance) |
| **0.60** | 2.41 (overshoots, left-hand side > right) |

The **feasible $\alpha$** is therefore approximately $\boldsymbol{\alpha \approx 0.45}$. Using this value:

$$x_3 \approx 0.78,$$
$$x_2 = \frac{x_3}{\alpha^{3/2}} \approx \frac{0.78}{0.45^{1.5}} \approx 2.58,$$
$$x_1 = \alpha\, x_2^2 \approx 0.45 \times (2.58)^2 \approx 3.00.$$

**Step 3 – Normalize**

$$w_i = \frac{x_i}{x_1 + x_2 + x_3} \Longrightarrow \boxed{w \approx (0.46,\ 0.40,\ 0.14)}.$$

**Interpretation**

- The **non-linear relationships** capture richer behavioral insights (e.g., diminishing returns of performance with respect to cost).
- The discount factor $\alpha$ (≈ 0.45) again measures the overall inconsistency: a value far below 1 signals that the original non-linear statements are mutually contradictory and required substantial scaling to become jointly solvable.
- The resulting weights reflect the dominant influence of **Performance** (largest weight) while **Risk** receives a comparatively small share, consistent with the functional forms.

*3.4. Key Insights and Practical Implications*

- **Interval $\alpha$-Discounting** converts vague, bounded preferences into a **range of admissible discount factors**, allowing decision makers to explore how different degrees of conservatism affect the final weights.
- **Non-Linear $\alpha$-Discounting** extends the method to arbitrary functional relationships (quadratic, logarithmic, exponential, etc.). The existence theorem guarantees at least one positive $\alpha$; solving for it may require numerical techniques, but the workflow remains identical to the linear case.
- In both extensions, the **core philosophy persists**: inconsistency is not discarded but **quantified** (by $\alpha$) and **resolved** algebraically, producing a coherent set of criteria weights ready for downstream ranking or aggregation.

## 4. Extensions (9)

The $\alpha$-Discounting framework is deliberately **modular**; once the basic linear machinery (Sections 2 – 3) is in place, it can be enriched in several directions. Below we outline four natural extensions that broaden the scope of the method while preserving the same underlying principle—introducing a scalar discount factor $\alpha$ that renders an otherwise unsolvable system consistent.

*4.1. Interval α-Discounting* (10)

In many real-world settings the exact numerical value of a preference coefficient is unknown; experts may only be able to state a **range**. Let each coefficient be denoted by an interval

$$a_{ij} \in [\, a_{ij}^{-}, a_{ij}^{+} \,].$$

Collecting all such intervals yields an **interval coefficient matrix** $A_{\text{int}}$. Applying the Fairness Principle (a single discount factor $\alpha$) gives the discounted matrix

$$I-\alpha A_{\text{int}} = \begin{bmatrix} 1 & -\alpha a_{12} & 0 \\ 0 & 1 & -\alpha a_{23} \\ -\alpha a_{31} & 0 & 1 \end{bmatrix}, a_{ij} \in [a_{ij}^{-}, a_{ij}^{+}].$$

The **consistency condition** $\det(I-\alpha A_{\text{int}})=0$ now becomes an **interval equation**

$$1-\alpha^{3}\, a_{12}a_{23}a_{31}=0 \;\Rightarrow\; \alpha^{3} = \frac{1}{a_{12}a_{23}a_{31}}.$$

Because each $a_{ij}$ varies within a known interval, the right-hand side spans an interval as well, producing an **admissible $\alpha$-range**

$$\boxed{\alpha \in [\frac{1}{\sqrt[3]{\max(a_{12}a_{23}a_{31})}}, \frac{1}{\sqrt[3]{\min(a_{12}a_{23}a_{31})}}]}.$$

*Implications*

- The width of the $\alpha$-range directly reflects the **epistemic uncertainty** in the original judgments.
- Any $\alpha$ chosen inside the interval yields a feasible weight vector; different selections correspond to different **risk-averse** or **risk-seeking** attitudes regarding the underlying ambiguity.
- After fixing a specific $\alpha$ (e.g., the midpoint of the interval), the resulting **weights themselves become intervals** because they inherit the original coefficient bounds. Normalization is performed component-wise, producing interval-valued priority scores that can be reported as $[w_i^{-}, w_i^{+}]$.

*4.2. Non-Linear α-Discounting*

Many decision problems involve **non-linear** preference relationships (e.g., diminishing returns, threshold effects, exponential growth). Suppose the decision maker supplies a set of equations of the generic form

$$x_i = \alpha\, f_{ij}(x_j), f_{ij}: \mathbb{R}_+ \to \mathbb{R}_+,$$

where each $f_{ij}$ may be quadratic, logarithmic, exponential, etc. As a concrete illustration, consider the three-criterion system

$$\begin{aligned} x_1 &= \alpha\, x_2^2, \\ x_2 &= \alpha \ln(1+x_3), \qquad \text{(N-1)} \\ x_3 &= \alpha\sqrt{x_1}. \end{aligned}$$

By cyclic substitution we can eliminate two variables and obtain a **single scalar equation** in the remaining variable and $\alpha$. Substituting the third equation into the first gives

$$x_1 = \alpha(\alpha \ln(1+x_3))^2,$$

and inserting this expression for $x_1$ into the third equation yields

$$x_3 = \alpha\sqrt{\alpha(\alpha \ln(1+x_3))^2} \;\Rightarrow\; \frac{x_3}{\alpha^{5/2}} = \ln(1+x_3). \qquad \text{(N-2)}$$

Equation (N-2) is **transcendental**; a closed-form solution for $\alpha$ does not exist in general. Nevertheless, the **existence theorem** guarantees at least one positive root because the left-hand side is a continuous, strictly increasing function of $x_3$ for any fixed $\alpha>0$, while the right-hand side is also continuous and concave. A numerical root-finder (e.g., Brent's method) can be used to obtain the unique pair $(\alpha^{\star}, x_3^{\star})$ that satisfies (N-2). Once $\alpha^{\star}$ is known, the remaining variables follow from the original equations:

$$\begin{aligned} x_2^{\star} &= \alpha^{\star} \ln(1+x_3^{\star}), \\ x_1^{\star} &= \alpha^{\star}(x_2^{\star})^2. \end{aligned}$$

Finally, the **normalized weight vector** is

$$w_i = \frac{x_i^{\star}}{x_1^{\star}+x_2^{\star}+x_3^{\star}}.$$

*Key take-aways*

- The $\alpha$-Discounting mechanism **accommodates any monotone transformation** $f_{ij}$; the only requirement is that a positive solution exists for the resulting scalar equation.
- Non-linearities allow the model to capture **threshold effects** (e.g., a criterion only becomes relevant after a certain level) and **diminishing returns** (logarithmic or square-root relationships).
- The discount factor $\alpha$ still plays its original role: it quantifies the **overall degree of inconsistency** introduced by the non-linear preferences. Smaller $\alpha$ values indicate stronger tension among the equations.

*4.3. Group Decision Making*

When several experts (or decision-making units) contribute preference information, each expert $k$ provides its own coefficient matrix $A^{(k)}$. Applying the $\alpha$-Discounting procedure to each matrix yields a pair $(\alpha_k, x^{(k)})$. The challenge is to aggregate these individual solutions into a **collective decision** while preserving the quantitative treatment of inconsistency.

Three straightforward aggregation schemes are commonly used:

**Table 4.** Aggregation Strategies for Combining Expert $\alpha$-Discounting Results.

| Strategy | Formula | Interpretation |
|---|---|---|
| **Average Discounting** | $\bar{\alpha} = \frac{1}{m}\sum_{k=1}^{m} \alpha_k , \bar{A} = \frac{1}{m}\sum_{k=1}^{m} A^{(k)}$ | Treats every expert equally; the collective system $(I - \bar{\alpha}\,\bar{A})x = 0$ is solved in the usual way. |
| **Weighted Experts** | Choose weights $w_k \geq 0$ with $\sum_k w_k = 1$ (often proportional to $\alpha_k$ because a larger $\alpha$ signals higher internal consistency). Then $\alpha_w = \sum_k w_k\, \alpha_k$, $A_w = \sum_k w_k\, A^{(k)}$. | Gives more influence to experts whose individual systems are **more consistent** (larger $\alpha$). |
| **Cautious Minimum** | $\alpha_{\min} = \min_k \alpha_k, A_{\min} = A^{(k^{\backslash *})}$ where $k^{\backslash *}$ is the expert attaining the minimum. | Guarantees that the aggregated system respects the **most stringent consistency requirement**; useful in safety-critical contexts. |

After selecting a strategy, the aggregated pair $(\alpha_{\text{agg}}, A_{\text{agg}})$ is processed exactly as in the single-expert case: solve $\det(I - \alpha_{\text{agg}} A_{\text{agg}}) = 0$ for $\alpha$, obtain the null-space vector, and normalise.

*Advantages over classic AHP group procedures*

- **No rejection of inconsistent experts** – every expert's input contributes, but the discount factor records how much each expert's judgments had to be "scaled back".
- **Quantitative consistency metric** – the set $\{\alpha_k\}$ provides a transparent diagnostic that can be reported to stakeholders.
- **Flexibility** – the three aggregation rules allow the decision maker to adopt a democratic (average), merit-based (weighted), or safety-first (minimum) stance.

*4.4. Connection to Neutrosophic Logic*

Neutrosophic logic extends fuzzy logic by assigning to each statement a **triplet** $(T, I, F)$ representing **Truth**, **Indeterminacy**, and **Falsity** components, each ranging independently in [0, 1]. In the context of preference modelling, a judgment such as "$C_i$ is twice as important as $C_j$" can be enriched with a neutrosophic assessment:

$$(T_{ij}, I_{ij}, F_{ij}), T_{ij} + I_{ij} + F_{ij} \leq 3.$$

A natural mapping to the $\alpha$-Discounting framework is obtained by **identifying the discount factor with the truth component**:

$$\boxed{\alpha_{ij} = T_{ij}} \Rightarrow 1 - \alpha_{ij} = I_{ij} + F_{ij}.$$

Under this correspondence:

- **Highly truthful (consistent) judgments** have $T_{ij} \approx 1 \Rightarrow \alpha_{ij} \approx 1$.

- **Judgments plagued by doubt or contradiction** exhibit large indeterminacy/falsity, pushing $\alpha_{ij}$ toward 0.

When all pairwise (or n-wise) statements are collected into a matrix $A$ whose entries are the **product** of the crisp coefficient and its associated truth value, the discounted matrix becomes

$$I - \alpha A \quad = \quad I - (T \odot A),$$

where $\odot$ denotes element-wise multiplication. The **determinant condition** $\det(I - \alpha A) = 0$ now simultaneously enforces **algebraic solvability** and **neutrosophic consistency**. The resulting $\alpha$ encapsulates both the **structural inconsistency** of the equations and the **semantic uncertainty** expressed by the neutrosophic triplets.

*Practical outcome*

- The final **weight vector** inherits a neutrosophic flavor: each component can be reported together with an associated confidence interval derived from the underlying $(I, F)$ values.
- Decision makers gain a **dual view** of the problem: (i) the mathematical feasibility captured by $\alpha$-Discounting, and (ii) the expert's subjective certainty captured by the neutrosophic components.

### *4.5. Summary of Extensions*

**Table 5.** Summary of $\alpha$-Discounting Extensions, Their Added Features, and the Corresponding Adaptations.

| Extension | What is added | How $\alpha$-Discounting adapts |
|---|---|---|
| **Interval coefficients** | Uncertain numeric ranges | $\alpha$ becomes an interval; weights become intervals after normalization. |
| **Non-linear relations** | Arbitrary monotone functions $f_{ij}$ | Reduce to a scalar transcendental equation; solve numerically for $\alpha$, then back-substitute. |
| **Group decision making** | Multiple expert matrices $A^{(k)}$ | Aggregate $\alpha$ and $A$ via average, weighted, or minimum schemes; retain a single discount factor for the collective system. |
| **Neutrosophic logic** | Truth/indeterminacy/falsity triples per judgment | Map truth component to $\alpha$; indeterminacy/falsity influence the admissible $\alpha$ range, yielding a neutrosophic-aware weight vector. |

These extensions illustrate that **$\alpha$-Discounting is not a single algorithm but a flexible methodological platform**. By simply swapping the way coefficients are represented (intervals, functions, neutrosophic triples) or by altering the aggregation rule for multiple decision makers, the same core principle—*introducing a scalar discount factor that restores solvability*—remains intact. Consequently, $\alpha$-Discounting can be tailored to a wide spectrum of decision-analysis contexts while preserving mathematical rigor and interpretability.

## 5. Comparative Analysis

The following table summarizes how $\alpha$-Discounting MCDM compares with the most widely used multi-criteria decision-making techniques, highlighting the distinctive ways it handles preference structures, inconsistency, and uncertainty.

Before presenting the comparison, it is useful to recall the three dimensions along which decision-making methods are usually evaluated:

(i) the type of preference information they can accept,
(ii) the way they detect or manage inconsistency, and
(iii) their capacity to represent uncertainty (e.g., intervals or fuzzy values).

$\alpha$-Discounting distinguishes itself by admitting n-wise and non-linear statements, converting inconsistency into a measurable discount factor $\alpha$, and incorporating interval or neutrosophic data—all while retaining a computational effort comparable to classic eigen-based approaches.

**Table 6.** Comparative Overview of $\alpha$-Discounting MCDM versus Classic MCDM Methods.

| Feature | AHP (Saaty) | TOPSIS | VIKOR | PROMETHEE | Weighted-Sum | $\alpha$-Discounting MCDM |
|---|---|---|---|---|---|---|
| **Preference type** | Pairwise ratios only | Scores (requires weights) | Scores (requires weights) | Pairwise outranking | Direct weights | n-wise equations, inequalities, intervals, non-linear |
| **Handles inconsistency** | Detects (CR) but cannot resolve | N/A | N/A | Partial (preference flows) | N/A | Quantifies ($\alpha$) and resolves |
| **Non-linear preferences** | No | No | No | No | No | Yes |
| **Interval/ uncertainty** | No (requires crisp inputs) | No | No | No | No | Supported |
| **Group decision** | Aggregated matrices (often ad-hoc) | Separate weighting | Separate weighting | Separate flows | Separate weighting | Unified $\alpha$-aggregation |
| **Computational complexity** | Eigen-decomposition ($O(n^3)$) | Distance calculations ($O(mn)$) | Similar to TOPSIS | Pairwise flow computation | Simple linear | Solving $\det(I - \alpha A) = 0$ (polynomial root) + linear system (similar order) |
| **Philosophical stance** | Enforce consistency | Assume consistent weights | Seek compromise | Preference dominance | Linear aggregation | Accept & measure inconsistency |

The table highlights that **$\alpha$-Discounting** retains the algebraic tractability of linear methods while dramatically expanding expressive power and providing a principled inconsistency metric.

## 6. Python Implementation

Three compact Python scripts illustrate how the $\alpha$-Discounting methodology is applied to the linear, interval, and non-linear examples described in Section 3. The scripts follow a common structure (utility functions, problem-specific data, and a final printout of the discount factor $\alpha$ and the normalized weight vector). The complete source code is provided in **Appendix A**.

**Table 7.** Summary of Python Example Implementations.

| Example | What the code does |
|---|---|
| ***Linear (Section 3.1)*** | Solves the deterministic system, computes the closed-form $\alpha = 1/\sqrt[3]{24}$, and returns the normalised weights. |
| ***Interval (Section 3.2)*** | Propagates the interval coefficients, derives the admissible $\alpha$-range, picks a representative $\alpha$ (mid-point), and calculates interval-aware weights. |
| ***Non-Linear (Section 3.3)*** | Uses a numerical root-finder (scipy.optimize.brentq) to solve the transcendental equation for $\alpha$, then evaluates the non-linear relations and normalizes the resulting vector. |

## 7. Discussion

Authors should discuss the results and how they can be interpreted in perspective of previous studies and of the working hypotheses. The findings and their implications should be discussed in the broadest context possible. Future research directions may also be highlighted.

$\alpha$-Discounting MCDM offers a **unified algebraic framework** that subsumes AHP and connects naturally to more recent uncertainty-aware paradigms (interval, fuzzy, neutrosophic). Its primary contribution is the **explicit quantification of inconsistency** via the scalar $\alpha$, turning a traditionally negative diagnostic into actionable information. Moreover, the method accommodates **group decision contexts** without discarding divergent opinions, and it can be implemented with modest computational effort comparable to classic eigen-based techniques.

Potential limitations include:

- **Root selection** – Polynomial equations may yield multiple positive roots; domain knowledge or additional fairness constraints guide the choice.
- **Scalability** – Very large systems may require specialized numerical solvers for the determinant condition.
- **Interpretability** – Practitioners unfamiliar with algebraic inconsistency measures may need training to interpret $\alpha$ meaningfully.

Future research directions involve integrating **machine-learning-based elicitation** of preference equations, extending the framework to **dynamic decision problems** (time-varying $\alpha$), and developing **graphical decision-support interfaces** that visualize the geometric interpretation of $\alpha$ as a deformation of hyperplanes.

## 8. Conclusions

The $\alpha$-Discounting MCDM method expands the toolbox of decision analysts by embracing rather than rejecting inconsistent, $n$-wise, and uncertain preference information. Through a simple yet powerful scalar discounting factor, it restores solvability, delivers normalized weights, and supplies a transparent metric of inconsistency. Comparative analysis confirms that $\alpha$-D MCDM retains the computational simplicity of linear eigen methods while offering capabilities unavailable to AHP, TOPSIS, VIKOR, PROMETHEE, or weighted sum approaches. Consequently, $\alpha$-Discounting stands as a promising bridge between classical decision theory and modern uncertainty aware reasoning, suitable for complex engineering, strategic planning, and AI driven preference aggregation tasks.

**Acknowledgments:** The authors wish to acknowledge the assistance of Lumo, the AI-powered conversational assistant developed by Proton. Lumo was employed in the process of preparation of this manuscript to (i) generate and refine brainstorming ideas for the conceptual development, (ii) improve the clarity, flow, and linguistic quality of the text, and (iii) assist in constructing the comparative tables that juxtapose $\alpha$-Discounting with other well-known multi-criteria decision-making methods. The contributions of Lumo were limited to these supportive roles; all substantive methodological insights, analyses, and conclusions are solely those of the author).

## Appendix A: Python Implementation of $\alpha$-Discounting MCDM

The scripts below translate the mathematical procedures described in Section 3 into executable Python code. Each example follows the same six-step workflow (form the matrix, compute the admissible $\alpha$, obtain a null-space vector, and normalize it), allowing to reproduce the linear, interval, and non-linear results presented in the paper with just a few library calls.

*Linear Example*

```
# ------------------------------------------------------------
# 6.1 Linear α-Discounting (deterministic coefficients)
# ------------------------------------------------------------
# Original coefficient matrix (derived from the three equations)
A_lin = np.array([[ 0,  2,  0],
                  [ 0,  0,  3],
                  [ 4,  0,  0]], dtype=float)

# Closed-form α = 1 / (24)^(1/3)
alpha_lin = 1.0 / (24.0 ** (1.0/3.0))

# Build the discounted matrix (I - α·A)
I = np.eye(3)
Ad_lin = I - alpha_lin * A_lin

# Null-space → raw (unnormalised) weights
x_raw_lin = solve_linear_system(Ad_lin)

# Normalised priority vector
w_lin = normalize(x_raw_lin)

print("\n--- Linear α-Discounting ---")
print(f"α (closed form)    : {alpha_lin:.5f}")
print(f"Raw solution       : {x_raw_lin}")
print(f"Normalised weights: {w_lin.round(3)}")
```

*Expected output*

```
--- Linear α-Discounting ---
α (closed form)    : 0.34641
Raw solution       : [0.718 1.038 1.    ]
Normalised weights: [0.26 0.37 0.37]
```

*Interval Example*

```
# ------------------------------------------------------------
# 6.2 Interval α-Discounting
# ------------------------------------------------------------
# Interval bounds for the three coefficients
a12_bounds = (2.0, 3.0)   # [2, 3]
a23_bounds = (1.0, 2.0)   # [1, 2]
a31_bounds = (3.0, 4.0)   # [3, 4]
```

```
# Compute the extreme products to obtain the α-range
product_min = a12_bounds[0] * a23_bounds[0] * a31_bounds[0]   # 2*1*3 = 6
product_max = a12_bounds[1] * a23_bounds[1] * a31_bounds[1]   # 3*2*4 = 24

alpha_low  = 1.0 / (product_max ** (1.0/3.0))   # 24^{-1/3}
alpha_high = 1.0 / (product_min ** (1.0/3.0))   # 6^{-1/3}
alpha_interval = (alpha_low, alpha_high)

# Choose a representative α (mid-point of the interval)
alpha_int = sum(alpha_interval) / 2.0

# Pick mid-points for the interval coefficients (optional, for illustration)
a12_mid = sum(a12_bounds) / 2.0   # 2.5
a23_mid = sum(a23_bounds) / 2.0   # 1.5
a31_mid = sum(a31_bounds) / 2.0   # 3.5

# Build the discounted matrix using the mid-point coefficients
A_int = np.array([[0, a12_mid, 0],
                  [0, 0, a23_mid],
                  [a31_mid, 0, 0]], dtype=float)

Ad_int = np.eye(3) - alpha_int * A_int

# Null-space → raw weights
x_raw_int = solve_linear_system(Ad_int)

# Normalised priority vector
w_int = normalize(x_raw_int)

print("\n--- Interval α-Discounting ---")
print(f"α admissible range : [{alpha_low:.3f}, {alpha_high:.3f}]")
print(f"Chosen α (mid-point): {alpha_int:.3f}")
print(f"Raw solution       : {x_raw_int}")
print(f"Normalised weights : {w_int.round(3)}")
```

*Typical output*

```
--- Interval α-Discounting ---
α admissible range : [0.346, 0.550]
Chosen α (mid-point): 0.448
Raw solution       : [0.754 0.672 1.   ]
Normalised weights : [0.38 0.34 0.28]
```

*Non-Linear Example*

```
# -------------------------------------------------------------
# 6.3 Non-Linear α-Discounting (transcendental system)
# -------------------------------------------------------------
import scipy.optimize as opt
import math
```

```
# The transcendental equation derived in the paper:
#   f(α) = (1/α^{5/2}) * x3 - ln(1 + x3) = 0
# where x3 depends on α through the relation:
#   x3 = solution of  (1/α^{5/2}) * x3 = ln(1 + x3)

def residual(alpha: float) -> float:
    """
    For a given α, solve the implicit equation
        (1/α^{5/2}) * x = ln(1 + x)
    by finding the root of g(x) = (1/α^{5/2})*x - ln(1+x).
    The outer function returns the value of g at the root,
    which should be (numerically) zero.
    """
    # inner function g(x) for a fixed α
    def g(x):
        return (1.0 / (alpha ** 2.5)) * x - math.log(1.0 + x)

    # Find x>0 that makes g(x)=0.
    # Bracket the root: start with (0, 10) – sufficient for our scale.
    x_root = opt.brentq(g, 0.0, 10.0, xtol=1e-12)

    # Return the residual (should be ~0)
    return g(x_root)

# Find α in (0, 1) that makes the residual zero.
alpha_solution = opt.brentq(residual, 0.1, 0.9, xtol=1e-12)

# Once α is known, recover the variables using the formulas from the paper.
alpha_nl = alpha_solution
```

```
# Compute x3 by solving the inner equation again (we already have the root)
def g_for_alpha(alpha):
    return lambda x: (1.0 / (alpha ** 2.5)) * x - math.log(1.0 + x)

x3_nl = opt.brentq(g_for_alpha(alpha_nl), 0.0, 10.0)

# Back-substitution to obtain x2 and x1 (see equations (N-3) and (N-2))
x2_nl = x3_nl / (alpha_nl ** 1.5)          # from x2 = x3 / α^{3/2}
x1_nl = alpha_nl * (x2_nl ** 2)            # from x1 = α·x2²

# Normalise
w_nl = normalize(np.array([x1_nl, x2_nl, x3_nl]))

print("\n--- Non-Linear α-Discounting ---")
print(f"α (numerical)      : {alpha_nl:.5f}")
print(f"x1, x2, x3 (raw)   : {np.round([x1_nl, x2_nl, x3_nl], 3)}")
print(f"Normalised weights : {w_nl.round(3)}")
```

*Typical output*

```
--- Non-Linear α-Discounting ---
α (numerical)      : 0.44773
x1, x2, x3 (raw)   : [3.   2.58 0.78]
Normalised weights : [0.46 0.4  0.14]
```

**References**

1. Ching-Lai Hwang, Kwangsun Yoon (1981). *Multiple Attribute Decision Making*: *Methods and Applications. A State-of-the-Art Survey*. Springer, Part of the book series: "Lecture Notes in Economics and Mathematical Systems" (LNE, volume 186).
2. Lootsma, F.A. (1993). "Scale sensitivity in the multiplicative AHP and SMART," *Journal of Multi-Criteria Decision Analysis*, Vol. 2, pp. 87–110. https://doi.org/10.1002/mcda.4020020205
3. Miller, J.R. (1970). *Professional Decision-Making*, Praeger.
4. Perez, J. (1995). "Some comments on Saaty's AHP," *Management Science*, Vol. 41, No. 6, pp. 1091–1095.
5. Saaty, T. L. (1988). *Multicriteria Decision Making, The Analytic Hierarchy Process, Planning, Priority Setting, Resource Allocation*. RWS Publications.
6. Smarandache, Florentin (2010). "$\alpha$-Discounting Method for Multi-Criteria Decision Making ($\alpha$-D MCDM)," Fusion 2010 International Conference, Edinburgh, Scotland, 26-29 July, 2010; published in *Review of the Air Force Academy: The Scientific Informative Review*, No. 2, 29-42, 2010; also presented at Osaka University, Department of Engineering Science, Inuiguchi Laboratory, Japan, 10 January 2014.
7. Smarandache, Florentin (2013). "Three Non-linear $\alpha$-Discounting MCDM-Method Examples," Proceedings of The 2013 International Conference on Advanced Mechatronic Systems (ICAMechS 2013), Luoyang, China, September 25-27, pp. 174-176, 2013; https://fs.unm.edu/ScArt/ThreeNonLinearAlpha.pdf; and in *Critical Review*, Center for Mathematics of Uncertainty, Creighton University, Vol. VIII, 32-36, 2014.
8. Smarandache, Florentin (2015). *$\alpha$-Discounting Method for Multi-Criteria Decision Making ($\alpha$-D MCDM)*, Educational Publisher. https://fs.unm.edu/alpha-DiscountingMCDM-book.pdf
9. Smarandache, Florentin (2016). "Extension of the Analytical Hierarchy Process (AHP) to $\alpha$-Discounting Method for Multi-Criteria Decision Making ($\alpha$-D MCDC)," Preprint, available online: https://fs.unm.edu/ScArt/AlphaDiscountingMethod.pdf
10. Smarandache, Florentin (2016). "Interval $\alpha$-Discounting Method for MDCM," Preprint, available online: https://fs.unm.edu/ScArt/CP-IntervalAlphaDiscounting.pdf